\documentstyle[12pt]{article}

\author{Roman R. Zapatrin}

\title{Algebraic duality for partially ordered sets}
\date{} 
\newcommand{\0}{{\bf 0}}
\newcommand{\1}{{\bf 1}}
\newcommand{\2}{{\bf 2}}
\newcommand{\poset}{{\cal POSET}}
\newcommand{\bcl}{{\cal BCL}}
\newcommand{\mor}{{\sf Mor}}
\newcommand{\lap}{\lambda_p}
\newcommand{\iup}{\upsilon_p}
\newcommand{\eqv}{\kern0.3em\Leftrightarrow\kern0.3em}
\newcommand{\implies}{\kern0.3em\Rightarrow\kern0.3em}

\newcommand{\proof}{\paragraph{Proof.}}

\newtheorem{lemma}{Lemma}

\begin{document}

\maketitle

\begin{center}
Department of Mathematics, SPb UEF, Griboyedova 30/32, \\
191023, St-Petersburg, Russia \\
{\em and}
\\ Division of Mathematics, Istituto per la Ricerca di Base, \\
I-86075, Monteroduni, Molise, Italy
\end{center}

\begin{abstract}
For an arbitrary partially ordered set $P$ its {\em dual} $P^*$ is
built as the collection of all monotone mappings $P\to\2$ where
$\2=\{0,1\}$ with $0<1$. The set of mappings $P^*$ is proved to be a
complete lattice with respect to the pointwise partial order. The
{\em second dual} $P^{**}$ is built as the collection of all morphisms
of complete lattices $P^*\to\2$ preserving universal bounds. Then it
is proved that the partially ordered sets $P$ and $P^{**}$ are
isomorphic.
\end{abstract}


\section*{Introduction}

The results presented in this paper can be considered as the
algebraic counterpart of the duality in the theory of linear
spaces. The outline of the construction looks as follows.

Several categories occur in the theory of partially ordered sets.
The most general is the category $\poset$ whose objects are
partially ordered sets and the morphisms are the monotone mappings.
Another category which will be used is
$\bcl$ whose objects are (bounded) complete lattices
and the morhisms are the lattice homomorphisms preserving
universal bounds.  Evidently $\bcl$ is the subcategory of $\poset$.

To introduce the algebraic duality (I use the term `algebraic' to
avoid confusion with the traditional duality based on order reversal)
the two element partially ordered set $\2$ is used:
\[ \2=\{0,1\}\quad,\quad 0<1 \]

Let $P$ be an object of $\poset$. Consider its dual $P^*$:
\begin{equation}\label{f129}
P^* = \mor_{\scriptscriptstyle \poset}(P,\2)
\end{equation}

The set $P^*$ has the pointwise partial order. Moreover, it is
always the complete lattice with respect to this partial order
(section \ref{s130}). Furthermore, starting from $P^* \in \bcl$
(bounded compete lattices) consider the set $P^{**}$ of all morphisms
in the appropriate category:
\begin{equation}\label{f129d}
P^{**} = \mor_{\scriptscriptstyle \bcl}(P^*,\2)
\end{equation}

And again, the set of mappings $P^{**}$ is pointwise partaially
ordered. Finally, it is proved in section \ref{s139} that $P^{**}$
is isomorphic to the initial partially ordered set $P$ (the
isomorphism lemma \ref{l141}):
\[ P^{**} \simeq P \]

The account of the results is organized as follows. First it is proved
that $P^*$ (\ref{f129}) is complete lattice. Then the
embeddings $p\to\lap$ and $p\to\iup$ of the poset $P$ into $P^*$
are built (\ref{f131d}). Then it is shown that the
principal ideals $[0,\lap]$ in $P^*$ are prime for all $p\in P$
(lemma \ref{l139}). Moreover, it is shown that there is
no more principal prime ideals in $P^*$. Finally, it is observed
that the principle prime ideals on $P^*$ are in 1-1 correspondence
with the elements of $P^{**}$.

\section{The structure of the dual space} \label{s130}

First define the pointwise partial order on the elements of $P^*$
(\ref{f129}). For any $x,y\in P^*$
\begin{equation}\label{f130p}
x\le y \eqv \forall p \in P \quad x(p)\le y(p)
\end{equation}

Evidently the following three statements are equivalent for $x,y\in
P^*$:
\begin{equation}\label{f130}
\left.
\begin{array}{lcl}
&x\le y & \cr
\forall p\in P \quad x(p)=1 &\Rightarrow & y(p)=1 \cr
\forall p\in P \quad y(p)=0 &\Rightarrow & x(p)=0
\end{array}
\right.
\end{equation}

To prove that $P^*$ is complete lattice, consider its arbitrary
subset $K\subseteq P^*$ and define the following mappings
$u,v:P \to \2$:
\begin{equation}\label{f131i}
u(p) = \left\lbrace \begin{array}{rcl}
1, & \exists k\in K & k(p)=1 \cr
0, & \forall k\in K & k(p)=0
\end{array}\right. \qquad
v(p) = \left\lbrace \begin{array}{rcl}
0, & \exists k\in K & k(p)=0 \cr
1, & \forall k\in K & k(p)=1
\end{array}\right.
\end{equation}

The direct calculations show that both $u$ and $v$ are monotone
mappings: $u,v\in P^*$ and
\[ u=\sup_{P^*}K \quad,\quad v=\inf_{P^*}K \]
which proves that $P^*$ is the complete lattice. Denote by \0,\1 the
universal bounds of the lattice $P^*$:
\[ \forall p\in P \quad \0(p)=0 \quad, \quad \1(p)=1 \]

Let $p$ be an element of $P$. Define the elements $\lap,\iup\in
P^*$ associated with $p$: for all $q\in P$
\begin{equation}\label{f131d}
\lap(q) = \left\lbrace \begin{array}{rcl}
0 &, & q\le p\cr
1 &, & \hbox{otherwise}
\end{array}\right. \qquad
\iup(q) = \left\lbrace \begin{array}{rcl}
1 &, & q\ge p \cr
0 &, & \hbox{otherwise}
\end{array}\right.
\end{equation}

\begin{lemma}\label{l132} For any $x\in P^*,\quad p\in\ P$
\begin{equation}\label{f132l}
\begin{array}{rcl}
x(p)=0 &\Leftrightarrow & x\le \lap \hbox{ \rm in } P^* \cr
x(p)=1 &\Leftrightarrow & x\ge \iup \hbox{ \rm in } P^*
\end{array}
\end{equation}
\end{lemma}

\proof Rewrite the left side of the first equivalency as \( \forall
q \quad q\le p \implies x(q)=0$, hence $\forall q \quad\lap(q)=0
\implies x(q)=0$, therefore $x\le \lap$ by virtue of (\ref{f130}).
The second equivalency is proved likewise. \hspace*{\fill}$\Box$\medskip

We shall focus on the 'inner' characterization of the elements
$\lap,\iup$ in mere terms of the lattice $P^*$ itself. To do it,
recall the necessary definitions.
\medskip

Let $L$ be a complete lattice. An element $a\in L$ is called {\sc
join-irreducible} ({\sc meet-irreducible}) if it can not be
represented as the join (resp., meet) of a collection of elements
of $L$ different from $a$. To make this definition more verifiable
introduce for every $a\in L$ the following elements of $L$:
\begin{equation}\label{f131}
\begin{array}{rcl}
\check{a} &=& \inf_L\{x\in L\mid\quad x>a\} \cr
\hat{a} &=& \sup_L\{y\in L\mid\quad y<a\}
\end{array}
\end{equation}
\noindent which do exist since $L$ is complete. Clearly,
$\check{a}\ge a \ge \hat{a}$ and the equivalencies
\begin{equation}\label{f132}
\begin{array}{rcl}
a\neq \check{a} &\Leftrightarrow & a\hbox{ is meet-irreducible} \cr
a\neq \hat{a} &\Leftrightarrow & a\hbox{ is join-irreducible}
\end{array}
\end{equation}
follow directly from the above definitions.

\begin{lemma} \label{l138} An element $w\in P^*$ is meet
irreducible if and only if it is equal to $\lap$ for some $p\in P$.
Dually, $v\in P^*$ is join irreducible iff $v=\iup$ for some $p\in
P$.
\end{lemma}

\proof First prove that every $\lap$ is meet irreducible. To do it
we shall use the criterion (\ref{f132}). Let $p\in P$. Define $u\in
P^*$ as:
\[ u(q) = \left\lbrace \begin{array}{rcl}
0 &,& q<p \cr
1 &,& \hbox{otherwise}
\end{array}
\right.
\]
then the following equivalency holds:
\begin{equation}\label{f135}
x\le y \eqv ( \forall q \quad x(q)=0 \implies q<p )
\end{equation}
Now let $x>\lap$, then $x(p)=1$ (otherwise (\ref{f132l}) would
enable $x\le\lap$). Then $x>\lap$ implies $x\ge\lap$, hence
$\forall q \quad x(q)=0 \implies q\le p$, although $q=p$ is excluded,
hence we get exactly the right side of (\ref{f135}). That means
that
\[ u=\inf_{P^*}\{x\mid x>\lap\} = {\check{\lambda}}_p \]
differs from $\lap$, hence $\lap$ is meet irreducible by virtue of
(\ref{f131}). The second dual statement is proved quite
analogously.

Conversely, suppose we have a meet irreducible $w\in P^*$, hence,
according to (\ref{f131}), there exists $p\in P^*$ such that
$\check{w}(p)\neq 0$ while $w(p)=0$. The latter means $w\le \lap$
for this $p$. To disprove $w<\lap$ rewrite $\check{w}(p)\neq 0$
as \( \lnot(\inf\{x\mid w<x\}\quad =\quad \lap) \)
which is equivalent to
\[ \exists y (\forall x \quad w<x \Rightarrow y\le x) \quad \&
\quad \lnot(y\le \lap) \]
In particular, it must hold for $x=\lap$, thus the assumption
$w<\lap$ implies \( \exists y \quad y\le \lap \kern0.4em \&
\kern0.4em \lnot(y\le \lap) \),
and the only remaining possibility for $w$ is to be equal to $\lap$.
\hspace*{\fill}$\Box$\medskip

\paragraph{Dual statement.} The join irreducibles of $P^*$ are the
elements $\iup, p\in P$ and only they.

\section{Second dual and the isopmorphism lemma} \label{s139}

Introduce the necessary definitions. Let $L$ be a lattice. An {\sc
ideal} in $L$ is a subset $K\subseteq L$ such that
\begin{itemize}
\item $k\in K, x\le k \implies x\in K$
\item $a,b\in K \implies a\lor b \in K$
\end{itemize}

Replacing $\le$ by $\ge$ and $\lor$ by $\land$ the notion of
{\sc filter} is introduced. An ideal (filter) $K\subseteq L$ is
called {\sc prime} if its set complement $L\setminus K$ is a
filter (resp., ideal) in $L$. Now return to the lattice $P^*$.

\begin{lemma}\label{l139} For any $p\in P$ both the principal ideal
$[0,\lap]$ and the principal filter $[\iup,1]$ are prime in $P^*$.
Moreover,
\[ [\iup,1] = P^*\setminus [0,\lap] \]
\end{lemma}

\proof Fix up $p\in P$, then for any $x\in P^*$ the value $x(p)$
is either 0 (hence $x\le \lap$) or 1 (and then $x\ge \iup$)
according to (\ref{f132l}). Since $\lap$ never equals
$\iup$ (because their values at $p$ are different), the sets
$[\iup,1]$ and $[0,\lap]$ are disjoint, which completes the proof.
\hspace*{\fill}$\Box$\medskip

The converse statement is formulated in the following lemma.

\begin{lemma} \label{l141} For any pair $u,v\in P^*$ such that
\begin{equation}\label{f141}
[0,u] = P^*\setminus [v,1]
\end{equation}
there exists an element $p\in P$ such that $u=\lap$ and $v=\iup$.
\end{lemma}

\proof It follows from (\ref{f141}) that $u$ and $v$ are not
comparable, therefore $u\land v < v$. Thus there exists $p\in P$
such that $(u\land v)(p) =0$ while $v(p) = 1$. Then (\ref{f132l})
implies $u\land v \le \lap$ and $v\ge \iup$. Suppose $v\neq \iup$,
then (\ref{f141}) implies $\iup\le u$, which together with $\iup\le
v$ implies $\iup \le u\land v \le \lap$ which never holds since
$\iup$ and $\lap$ are not comparable. So, we have to conclude that
$v=\iup$, thus $u=\lap$. \hspace*{\fill}$\Box$\medskip

Now introduce the {\sc second dual} $P^{**}$ as the set of all
homomorphisms of complete lattices $P^*\to \2$ preserving universal
bounds, that is, for any ${\bf p} \in P^{**}, K\subseteq P^*$
\[ \begin{array}{l}
{\bf p}(\sup_K) = \sup_{k\in K}{\bf p}(k) \cr
{\bf p}(\inf_K) = \inf_{k\in K}{\bf p}(k) \cr
{\bf p}(0) = 0\hbox{ ; }{\bf p}(1) = 1
\end{array} \]
with the pointwise partial order as in (\ref{f130p}).

Now we are ready to prove the following {\em isomorphism lemma}.

\begin{lemma} The partially ordered sets $P$ and $P^{**}$ are
isomorphic.
\end{lemma}

\proof Define the mapping $F:P\to P^{**}$ by putting
\[ F(p) = {\bf p}:\quad {\bf p}(x) = x(p) \qquad \forall x\in P^{**} \]

Evidently $F$ is the order preserving injection. To build the
inverse mapping $G:P^{**}\to P$, for any ${\bf p}\in P^{**}$ consider
the ideal ${\bf p}^{-1}(0)$ and the filter ${\bf p}^{-1}(1)$ in $P^*$
both being prime (see \cite{lt}, II.4). Let $u=\sup{\bf p}^{-1}(0)$
and $v=\inf{\bf p}^{-1}(1)$. Since ${\bf p}$ is the homomorphism of
complete lattices, $u\in {\bf p}^{-1}(0)$ and $v\in {\bf p}^{-1}(1)$,
hence ${\bf p}^{-1}(0) = [0,u]$ and ${\bf p}^{-1}(1) = [v,1]$. Applying
lemma \ref{l141} we see that there exists $p\in P$ such that $u=\lap$
and $v=\iup$. Put $G({\bf p}) = p$. The mapping $G$ is order preserving
and injective (since the different principal ideals have different
suprema). It remains to prove that $F,G$ are mutually inverse.

Let $p\in P$, consider $G(F(p))$. Denote ${\bf p} = F(p)$, then
${\bf p}^{-1}(1) = \{x\in P^*\mid x(p)=0\} = \{x\mid x\le \lap\}$. Thus
$\sup{\bf p}^{-1}(0) = \lap$, then $G\circ F = {\rm id}_P$ which
completes the proof. \hspace*{\fill}$\Box$\medskip

\section*{Concluding remarks}

     The results presented in this paper show that besides the well
known duality in partially ordered sets based on order reversal, we
can establish quite another kind of duality {\em \`a la} linear algebra.
As in the theory of linear topological spaces, we see that the
`reflexivity' expressed as $P=P^{**}$ can be achieved by appropriate
{\em definition} of dual space.

     We see that a general partially ordered set have the dual
space being a complete lattice. We also see that not every complete
lattice can play the r\^ole of dual for a poset. These complete
lattices can be characterized in terms of spaces with two closure
operations \cite{cr}. For the category of {\em ortho}posets this
construction was introduced in \cite{mayet}. Another approach to
dual spaces when they are treated as sets of two-valued measures
(in terms of this paper, as sub-posets of $P^{**}$) is in
\cite{tkadlec}. The main feature of the techniques suggested in the
present paper is that all the constructions are formulated in mere
terms of partially ordered sets and lattices.

The work was supported by the RFFI research grant (97-14.3-62). 
The author acknowledges the financial support from the Soros 
foundation (grant A97-996) and the research grant "Universities 
of Russia".

\end{document}